\documentclass{article}
\usepackage{amssymb,amsmath,mathrsfs,latexsym,amscd}
\usepackage{exscale,latexsym,amsthm,graphics}

\usepackage{makeidx}
\makeindex
\usepackage{epsfig}
\usepackage{hyperref}
\usepackage{txfonts}
\newtheorem{thm}{Theorem}[section]
\newtheorem{cor}[thm]{Corollary}

\newtheorem{lemma}[thm]{Lemma}
\newtheorem{prop}[thm]{Proposition}
\newtheorem{proposition}[thm]{Proposition}

\newtheorem{remark}[thm]{Remark}

\def\bM {{\mathbb M}}

\def\N {{\mathbb N}}

\def\qed{{\hfill $\Box$ \bigskip}}
\def\N {{\mathbb N}}
\def\R {{\mathbb R}}

\def\EE{{\mathbb E}}
\def\P{{\mathbb P}}
\newcommand{\F}{\mathcal{F}}

\def\E{{\mathcal E}}
\def\P{{\mathbb P}}

\numberwithin{equation}{section}
\begin{document}

\noindent
{{\Large\bf Non-symmetric distorted Brownian motion: strong solutions, strong Feller property and 
non-explosion results}\footnote{This research was supported by DFG through Grant Ro 1195/10-1 and by NRF-DFG Collaborative Research program and Basic Science Research Program 
through the National Research Foundation of Korea (NRF-2012K2A5A6047864 and NRF-2012R1A1A2006987).}}

\bigskip
\noindent
{\bf Michael R\"ockner},
{\bf Jiyong Shin},
{\bf Gerald Trutnau}
\\

\noindent
{\small{\bf Abstract.} Using elliptic regularity results in weighted spaces, stochastic calculus and the theory of non-symmetric Dirichlet forms, we first show 
weak existence of non-symmetric distorted Brownian motion for any starting point in some domain $E$ of $\R^d$, where $E$ is explicitly given as the points of strict 
positivity of the unique continuous version of the density to its invariant measure. This non-symmetric distorted Brownian motion is also proved to be strong Feller. Non-symmetric distorted Brownian motion is a singular diffusion, i.e. a diffusion that typically has 
an unbounded and discontinuous drift. Once having shown weak existence, we obtain from a result of \cite{KR} that the constructed 
weak solution is indeed strong and weakly as well as pathwise unique up to its explosion time. As a consequence of our approach, we can use the theory of Dirichlet forms 
to prove further properties of the solutions. For example, we obtain new non-explosion criteria for them. 
We finally present concrete existence and non-explosion results for non-symmetric distorted Brownian motion related 
to a class of Muckenhoupt weights and corresponding divergence free perturbations.\\ 



\section{Introduction}
In this paper we are concerned with the non-symmetric Dirichlet form given by (the closure of)
\begin{equation}\label{eq0.1}
 \mathcal{E}(f,g):=\frac{1}{2} \int_{\R^d} \langle \nabla f, \nabla g \rangle \;\mathrm{d} m - \int_{\R^d}\langle B, \nabla f \rangle g \;\mathrm{d}m\;, \quad f,g \in C_0^\infty(\R^d)\;,
\end{equation}
on $L^2(\R^d,m)$, $m:=\rho\; \mathrm{d}x$, and the corresponding stochastic differential equation (SDE)
\begin{equation}\label{eq0.2}
 X_t=x+W_t+\int_0^t \left( \frac{\nabla \rho}{2 \rho}+B\right)\left(X_s\right)\mathrm{d} s\;, \quad t<\zeta\;,
\end{equation}
where $x \in \R^d$, $\zeta$ is the lifetime (=explosion time).
Our conditions on $\rho$ and $B$ are formulated as Hypotheses (H1)-(H3) in Section \ref{S2} below. 

It is well-known that starting with \eqref{eq0.1} by Dirichlet form theory one can construct a weak solution to \eqref{eq0.2} for quasi-every starting point $x\in \R^d$, and usually there is no analytic characterization (in terms of $\rho$ and $B$)
of the set of \lq\lq allowed\rq\rq\ starting points. \\
In case $B\equiv0$, it was however shown in \cite{AKR} (see also \cite{B},\cite{FG1},
for extensions of this result to other situations), 
that \eqref{eq0.2} has a weak solution for every $x\in\{\tilde{\rho}>0\}$ in the sense of the martingale problem, where $\tilde{\rho}$
is the continuous version of $\rho$ (which exists as a consequence of (H1)) and that for such starting points the process $X_t$ stays in $\{\tilde{\rho}>0\}$ before its lifetime $\zeta$. 
The identification of (\ref{eq0.2}) with $B\equiv0$ for any $x\in\{\tilde{\rho}>0\}$ in the sense of a weak solution of an SDE related to the form in (\ref{eq0.1}) has been worked out as a part of a general framework in \cite[Section 4]{ShTr13a}.

The first aim of this paper is to 
generalize these results to $B\nequiv0$, i.e. to the non-symmetric case (see Remark \ref{r2.2}). The proof follows ideas from \cite{AKR}, i.e., in particular, we first construct a {\it strong Feller} semigroup of kernels on $\{\tilde{\rho}>0\}$, which are versions of the operator semigroup $(T_t)_{t>0}$ associated to the closure of (\ref{eq0.1})  and to be the transition semigroup of the corresponding process. However, compared to \cite{AKR} a number of modifications of the arguments 
there are required. For example, one observation is that the elliptic regularity results in weighted spaces from \cite{AKR} extend to the non-symmetric case.
The corresponding result is formulated 
as Theorem \ref{t3.6} in Section \ref{S3} below.

It is well-known by \cite[Theorem 2.1]{KR} (see also \cite{FF}, \cite{Zh}) that for every $x \in \{\tilde{\rho}>0\}$ there exists a strong solution (i.e. adapted to the filtration generated by $(W_t)_{t\geq0}$) to \eqref{eq0.2}, which is 
pathwise and weak unique. Hence this solution coincides with our weak solution (which is hence a strong solution) from Theorem \ref{t3.6}. Thus we have identified the Dirichlet form associated to the (strong Feller) Markov processes, given by the 
laws $\mathbb{P}_x, x\in \{\tilde{\rho}>0\}$, of these strong solutions, to be the closure of \eqref{eq0.1}. 
We emphasize that such an identification is generically non-trivial, since starting from the Markov process 
one can usually only identify the generator of the associated Dirichlet form on nice functions as those in $C_0^2(\R^d)$. 
But there are in general many Dirichlet forms extending (\ref{eq0.1}) with this property which are different from the 
one obtained as the closure of (\ref{eq0.1}). In general, if only one such extension exists, which is then necessarily 
the closure, one says \lq\lq Markov uniqueness holds\rq\rq. But the latter is unknown in our case. 
As a consequence of the aforementioned identification, we can apply the theory of Dirichlet forms to obtain further properties of the solutions to \eqref{eq0.2} for every
starting point in $\{\tilde{\rho}>0\}$.

In this paper, as our second aim, we concentrate on proving non-explosion results for \eqref{eq0.2} using Dirichlet form theory, 
which means (cf. Remark \ref{conserve}) that the process started in $x \in \{\tilde{\rho}>0\}$ will neither go to infinity
nor hit any point in $\{\tilde{\rho}=0\}$ in finite time. Non-explosion criteria from Dirichlet form theory are of analytic nature and different from the usual ones known from the theory of SDE (e.g. the one proved in \cite{KR},
see Remark \ref{r4.2} (ii) below), but very useful in applications.

Finally, we present a number of concrete applications where the density $\rho\left(=\frac{\mathrm{d}m}{\mathrm{d}x}\right)$ is in certain Muckenhoupt classes. Our main result here is Theorem \ref{Muckenhouptmain}.

The organization of this paper is as follows. After this introduction in Section \ref{S2}, we recall some important elliptic regularity results for the Kolmogorov operator corresponding to \eqref{eq0.2}, i.e. the generator of the
Dirichlet form \eqref{eq0.1}, under the assumption (H1) on $\rho$ and (H2) on $B$. Subsequently, we present their analytic consequences, in particular, obtain the strong Feller semigroup of kernels (transition semigroup) mentioned above (see Proposition \ref{prop1.5}). In Section \ref{S3} we construct the weak solutions
of \eqref{eq0.2} for every $x\in\{\tilde{\rho}>0\}$. In Section \ref{S4} we show that by \cite[Theorem 2.1]{KR} these solutions are strong,  pathwise and weak unique. 
Section \ref{S5} is devoted to the mentioned applications.

\section{Elliptic regularity and construction of a diffusion process}\label{S2}
For $E \subset \R^d$ open with Borel  $\sigma$-algebra $\mathcal{B}(E)$, we denote the set of all $\mathcal{B}(E)$-measurable $f : E \rightarrow \R$ which are bounded, or nonnegative by $\mathcal{B}_b(E)$, $\mathcal{B}^{+}(E)$ respectively. $L^q(E, \mu)$, $q \in[1,\infty]$ are the usual $L^q$-spaces equipped with $L^{q}$-norm $\| \cdot \|_{L^q(E,\mu)}$ with respect to the  measure $\mu$ on $E$, $\mathcal{A}_b$ : = $\mathcal{A} \cap \mathcal{B}_b(E)$ for $\mathcal{A} \subset L^q(E,\mu)$,  and $L^{q}_{loc}(E,\mu) := \{ f \,|\; f \cdot 1_{U} \in L^q(E, \mu),\,\forall U \subset E, U \text{ relatively compact open} \}$, where $1_A$ denotes the indicator function of a set $A$. Let $\nabla f : = ( \partial_{1} f, \dots , \partial_{d} f )$  and  $\Delta f : = \sum_{j=1}^{d} \partial_{jj} f$ where $\partial_j f$ is the $j$-th weak partial derivative of $f$ and $\partial_{jj} f := \partial_{j}(\partial_{j} f) $, $j=1, \dots, d$.  As usual $dx$ denotes Lebesgue measure on $\R^d$ and the Sobolev space $H^{1,q}(E, dx)$, $q \ge 1$ is defined to be the set of all functions $f \in L^{q}(E, dx)$ such that $\partial_{j} f \in L^{q}(E, dx)$, $j=1, \dots, d$, and 
$H^{1,q}_{loc}(E, dx) : =  \{ f  \,|\;  f \cdot \varphi \in H^{1,q}(E, dx),\,\forall \varphi \in  C_0^{\infty}(E)\}$. 
Here $C_0^{\infty}(E)$ denotes the set of all infinitely differentiable functions with compact support in $E$. We also denote the set of continuous functions on $E$, the set of continuous bounded functions on $E$, the set of compactly supported continuous functions in $E$ by $C(E)$, $C_b(E)$, $C_0(E)$, respectively. $C_{\infty}(E)$ denotes the space of continuous functions on $E$ which vanish at infinity. 
We equip $\R^d$ with the Euclidean norm $\| \cdot \|$ with corresponding inner product $\langle \cdot, \cdot \rangle$ and write $B_{r}(x): = \{ y \in \R^d \ | \ \|x-y\| < r  \}$, $x \in \R^d$. The closure of $B \subset \R^d$ is denoted by $\overline{B}$.\\   

We shall assume (H1)-(H3) below throughout up to including section 3:
\begin{itemize}
\item[(H1)]
$\rho = \xi^2$, $\xi \in H^{1,2}_{loc}(\R^d, dx)$,  $\rho > 0 \ \ dx$-a.e. and 
\[
\frac{\| \nabla \rho \|}{\rho} \in L^{p}_{loc} (\R^d, m), \quad m := \rho dx,
\]
$p:=(d + \varepsilon) \vee 2$ for some $\varepsilon >0$. 
\end{itemize}
By (H1) the symmetric positive definite bilinear form
\[
\E^0(f,g) : = \frac{1}{2} \int_{\R^d} \langle \nabla f ,  \nabla g \rangle \ dm , \quad f, g \in C_0^{\infty}(\R^d)
\]
is closable in $L^2(\R^d,m)$ and its closure $(\E^0,D(\E^0))$ is a symmetric, strongly local, regular Dirichlet form. We further assume
\begin{itemize}
\item[(H2)] $B : \R^d \to \R^d, \ \|B\| \in L_{loc}^{p}(\R^d, m)$ where $p$ is the same as in (H1) and
\begin{equation}\label{eq1.1}
\int_{\R^d} \langle B,\nabla f \rangle \ dm = 0, \quad \forall f \in C_0^{\infty}(\R^d),
\end{equation} 
\end{itemize}
and
\begin{itemize}
\item[(H3)]
\[
\left| \int_{\R^d} \langle B,\nabla f \rangle \ g \ \rho \ dx \right| \le c_0 \ \E^0_1 (f,f)^{1/2} \ \E^0_1 (g,g)^{1/2}, \quad \forall f,g \in C_0^{\infty}(\R^d),
\]
where $c_0$ is some constant (independent of $f$ and $g$) and $\E_{\alpha}^0(\cdot,\cdot):=\E^0(\cdot,\cdot) + \alpha (\cdot,\cdot)_{L^2(\R^d,m)}$, $\alpha > 0$.
\end{itemize}
Next, we consider the non-symmetric bilinear form 
\begin{equation}\label{df}
\E(f,g) : = \frac{1}{2} \int_{\R^d} \langle \nabla f , \nabla g \rangle \ dm - \int_{\R^d} \langle B,\nabla f \rangle \ g \  dm, \quad
 f, g \in C_0^{\infty}(\R^d)
\end{equation}
in $L^2(\R^d,m)$. 
Then by (H1)-(H3) $(\E, C_0^{\infty}(\R^d))$ is closable in $L^2(\R^d,m)$ and the closure $(\E,D(\E))$ is a non-symmetric Dirichlet form (cf. \cite[II. 2.d)]{MR}). Let  $(T_t)_{t > 0}$ (resp. $(\hat{T}_t)_{t > 0}$) and $(G_{\alpha})_{\alpha > 0}$ (resp. $(\hat{G}_{\alpha})_{\alpha > 0}$ ) be the $L^2(\R^d, m)$-semigroup (resp. cosemigroup) and resolvent (resp. coresolvent) associated to $(\E,D(\E))$ and $(L,D(L))$ (resp. $(\hat{L},D(\hat{L}))$) be the corresponding generator (resp. cogenerator) (see \cite[Diagram 3, p. 39]{MR}). Using properties (H2) and \cite[I. Proposition 4.7]{MR} (cf. also \cite[II 2.d)]{MR}), it is straightforward to see that $(T_t)_{t>0}$ as well as $(\hat{T}_t)_{t>0}$ are submarkovian. Here an operator $S$ is called submarkovian if $0 \le f \le 1$ implies $0 \le Sf \le 1$.
It is then further easy to see that $(T_t)_{t >0}$ (resp. $(G_{\lambda})_{\lambda > 0}$) restricted to $L^1(\R^d,m) \cap L^{\infty}(\R^d,m)$ can be extended to strongly continuous contraction semigroups (resp. strongly continuous contraction resolvents) on all $L^r(\R^d,m)$, $r \in [1,\infty)$ (see \cite[I.1]{MR} for the definition of a strongly continuous contraction semigroup (resp. resolvent)). We denote the corresponding operator families again by $(T_t)_{t > 0}$ and $(G_{\lambda})_{\lambda > 0}$ and let $(L_r, D(L_r))$ be the corresponding generator on $L^r(\R^d,m)$.
Since by (H1), (H2), $\left\| \frac{\nabla \rho}{2 \rho} \right\|$, $\|B\| \in L^{p}_{loc}(\R^d,m)$, we get  $C_0^{\infty}(\R^d) \subset D(L_r)$ for any $r \in [1,p]$ and
\begin{equation}\label{eq1.3}
L_r u = \frac{1}{2}\Delta u +  \langle \frac{\nabla \rho}{2 \rho}+B, \nabla u \rangle, \quad u \in C_0^{\infty}(\R^d), \quad r \in [1,p].
\end{equation}

Let us first state an elliptic regularity result (cf. \cite[Theorem 1 (iii)(b)]{BKR1}, \cite[Remark 2.15]{BKR2}). Its consequences in the symmetric case were discussed in \cite{AKR}. Likewise the Corollaries \ref{cor1.2}, \ref{cor1.3}, \ref{cor1.4}, and Remark \ref{rem1.5} below can be obtained.

\begin{prop}\label{prop1.1} 
Let $E$ be an open set in $\R^d$ and $A : E \to \R^d$, $c : E \to \R$ Borel measurable maps. Suppose $\mu$ is a (signed) Radon measure on $E$ and $f \in L_{loc}^1 (E, dx)$ such that $\| A \|, \ c \in L^1_{loc}(E,\mu)$ and
\[
\int Nu(x) \ \mu(dx) = \int u(x) \ f(x) \ dx, \quad \forall u \in C_0^{\infty} (E),
\]
where
\[
N u(x) : = \Delta u(x) + \langle A(x), \nabla u(x) \rangle + c(x) \ u(x).
\]
If for some $\tilde p > d$, $\|A\| \in L^{\tilde p}_{loc}(E, \mu)$, $c \in L^{{\tilde p}d/({\tilde p}+d)}_{loc}(E, \mu)$, and $f \in L^{{\tilde p}d/({\tilde p}+d)}_{loc}(E, dx)$, then $\mu = \rho dx$ with $\rho$ continuous and 
\[
\rho \in H^{1,{\tilde p}}_{loc} (E, dx) \ \Big( \subset C^{1-d/{\tilde p}}_{loc}(E) \Big),
\]
where $C^{1-d/{\tilde p}}_{loc}(E)$ denotes the set of all locally H\"{o}lder continuous functions of order $1-d/{\tilde p}$ on $E$. If $E_0 : = E \cap \{\rho >0\}$ and moreover $f,c \in L^{\tilde p}_{loc}(E_0)$, then for any open balls $B'$, $B$ with $B'\subset \overline{B'}\subset B \subset \overline{B} \subset E_0$ there exists $c_B \in (0, \infty)$ (independent of $\rho$ and $f$) such that 
\[
\|\rho\|_{H^{1,{\tilde p}}(B',dx)} \le c_B \left(\|\rho\|_{L^1(B,dx)} + \|f\|_{L^{\tilde p}(B,dx)} \right).
\]
\end{prop}
\begin{remark}\label{r2.2}
 At first sight the assumption that the drift in \eqref{eq0.2} or the first order coefficient in \eqref{df} is of type $b:=\frac{\nabla \rho}{2 \rho}+B$ looks rather special. But the $L^p_{loc}(\R^d,m)$ condition makes it very natural, because the special form of $b$ follows, if one considers the operator
 \begin{equation*}
  Lu:=\frac12\Delta u + \langle b, \nabla u \rangle\;, \quad u \in C^\infty_0(\R^d)\;,
 \end{equation*}
and assumes that it has an infinitesimally (not necessarily probability) invariant measure $m$, i.e. $m$ is a nonnegative Radon measure $m$ on $\R^d$, such that $b \in L^p_{loc}(\R^d,m)$ and
\begin{equation*}
 \int Lu \ \mathrm{d} m =0, \quad \forall u \in C^\infty_0(\R^d),
\end{equation*}
since then it follows by Proposition \ref{prop1.1} that $m= \rho \mathrm{d}x$ and that $\rho$ satisfies (H1).\\
Defining 
\begin{equation*}
 B:=b - \frac{\nabla \rho}{2 \rho}\;,
\end{equation*}
it satisfies (H2). So, we have the above decomposition in a natural way. 
\end{remark}

\begin{cor}\label{cor1.2} 
$\rho$ is in $H^{1,p}_{loc}(\R^d,dx)$ and $\rho$ has a continuous $dx$-version in $C_{loc}^{1-d/p}(\R^d)$.
\end{cor}
\proof
By \eqref{eq1.1}, \eqref{eq1.3} and integration by parts, we obtain
\[
\int Lu \ \mathrm{d}m = 0, \quad \forall u \in C_0^{\infty}(\R^d).
\]
Since $\frac{\| \nabla \rho \|}{\rho}, \|B\| \in L_{loc}^{p}(\R^d, m)$, the assertion follows by Proposition \ref{prop1.1} applied with $\tilde p =p$.
\qed

From now on, we shall always consider the continuous $dx$-version of $\rho$  and denote it also by $\rho$.
\begin{cor}\label{cor1.3}
Let $\lambda > 0$.
Suppose $g \in L^r(\R^d, m)$, $r \in [p,\infty)$. Then
\[
\rho \ G_{\lambda} g \in H_{loc}^{1,p}(\R^d,dx)
\]
and for any open balls $B'\subset \overline{B'}\subset B \subset \overline{B}  \subset \{\rho > 0\}$ there exists $c_{B,\lambda} \in (0,\infty)$, independent of $g$, such that 
\begin{equation}\label{conresol}
\| \ \rho \ G_{\lambda} g \ \|_{H^{1,p}(B',dx)} \le c_{B,\lambda} \ \Big(\|G_{\lambda} g \|_{L^1(B,m)} + \| g \|_{L^{p}(B,m)} \Big).
\end{equation}
\end{cor}
\proof
Let $g \in C_0^{\infty}(\R^d)$. Then we have
\[
\int (\lambda - \hat{L})u \ G_{\lambda} g \ \rho \ dx = \int u \ g \ \rho \ dx, \quad \forall u \in C_0^{\infty}(\R^d),
\]
where
\[
\hat{L} u = \frac{1}{2} \Delta u +  \langle \frac{\nabla \rho}{2 \rho} - B, \nabla u \rangle.
\]
Now we apply Proposition \ref{prop1.1} with $\mu = - \frac{1}{2} \rho G_{\lambda} g dx$ and $N = -2 (\lambda - \hat{L})$ and $f=g \rho$ to prove the assertion for $g \in C_0^{\infty}(\R^d)$. Since $C_0^{\infty}(\R^d)$ is dense in ($L^r(\R^d, m)$, $\| \cdot \|_{L^r(\R^d, m)})$, $r \in [1,\infty) $, the assertion for general $g \in L^r(\R^d, m)$ follows by continuity and \eqref{conresol}.
\qed

\begin{remark}
By \cite[I. Corollary 2.21]{MR}, it holds that $(T_t)_{t>0}$ is analytic on $L^2(\R^d,m)$. By Stein interpolation (cf. e.g. \cite[Lecture 10, Theorem 10.8]{AV}) $(T_t)_{t>0}$ is also analytic on $L^r(\R^d,m)$ for all $r \in (2, \infty)$. We would like to thank Hendrik Vogt for pointing this out to us as well as a misprint in the mentioned Theorem 10.8. There $\theta_{\tau}$ should be defined as $\tau \cdot \theta$ and not as $(1- \tau)\cdot \theta$.
\end{remark}

\begin{cor}\label{cor1.4}
Let $t>0$, $r \in [p,\infty)$. 
\begin{itemize}
\item[(i)] Let $u \in D(L_r)$. Then 
\[
\rho \ T_t u \in H^{1,p}_{loc} (\R^d, dx)
\]
and for any open balls $B'\subset \overline{B'}\subset B \subset \overline{B}  \subset \{ \rho > 0\}$ there exists $c_B \in (0,\infty)$ (independent of $u$ and $t$) such that
\begin{eqnarray}\label{eq1.4}
\| \rho \ T_t u\|_{H^{1,p}(B',dx)} &\le& c_B \left( \| T_t u \|_{L^1(B,m)} + \| T_t (1-L_r) u  \|_{L^p(B,m)}   \notag    \right)\\
&\le& c_B \left( m(B)^{\frac{r-1}{r}} \| u \|_{L^r (\R^d,m) }  + m(B)^{\frac{r-p}{rp}} \| (1-L_r) u \|_{L^r(\R^d,m)}    \right).
\end{eqnarray}

\item[(ii)] Let $f \in L^r(\R^d,m)$. Then the above statements still hold with \eqref{eq1.4} replaced by
\[
\|\rho \ T_t f\|_{H^{1,p}(B',dx)} \le \tilde{c}_B \ (1+ t^{-1}) \| f \|_{L^r(\R^d,m)},
\]
where $\tilde{c}_B \in (0,\infty)$ (independent of $f$, $t$).
\end{itemize}
\end{cor}

\begin{remark}\label{rem1.5}
By \eqref{eq1.4} and Sobolev imbedding, for $r \in [p, \infty)$, $R>0$ the set
\[
\{T_t u  \ | \ t>0, \ u \in D(L_r), \ \|u\|_{L^r(\R^d,m)} + \|L_r u\|_{L^r(\R^d, m)} \le R \} 
\]
is equicontinuous on $\{\rho>0\}$.
\end{remark}

From now on, we shall keep the notation
\[
E : = \{\rho > 0 \}.
\] 
By Corollaries \ref{cor1.2}, \ref{cor1.3}, \ref{cor1.4} and Remark \ref{rem1.5}, exactly as in \cite[section 3]{AKR}, we obtain the existence of a transition kernel density $p_t(\cdot,\cdot)$ on the open set $E$
such that
\[
P_t f(x) : = \int_E f(y) p_t(x,y) \ m(dy), \quad x \in E, \ t > 0
\]
is a (temporally homogeneous) submarkovian transition function (cf. \cite[1.2]{CW}) and an $m$-version of $T_t f $ for any $f \in \cup_{r \ge p} L^r(E,m)$. Moreover, letting $P_0 : = id$, it holds
\begin{equation}\label{eq2.5}
P_t f \in C(E) \quad \forall f \in \cup_{r \ge p} L^r(E,m)
\end{equation}
and
\begin{equation}\label{eq2.6}
\lim_{t \to 0} P_{t+s} f(x) = P_s f(x) \quad \forall s \ge 0, \ x \in E, \ f \in C_0^{\infty}(\R^d).
\end{equation}
By a 3$\varepsilon$-argument \eqref{eq2.6} extends to $C_0(\R^d)$.
Similarly, since for $\lambda > 0$, $f \in L^{p}( E, m)$, $G_{\lambda} f$ has a unique continuous $m$-version on $E$ by Corollary \ref{cor1.3} as in \cite[Lemma 3.4, Proposition 3.5]{AKR}, we can find $(R_{\lambda})_{\lambda > 0}$ with resolvent kernel density $r_{\lambda}(\cdot, \cdot)$ defined on $E \times E$ such that
\[
R_{\lambda} f(x) : = \int f(y)\, r_{\lambda} (x,y) \ m(dy),  \quad x \in E, \ \lambda > 0,
\]
satisfies 
\begin{equation}\label{strongresol}
R_{\lambda}f  \in C(E) \ \text{and} \ R_{\lambda} f = G_{\lambda} f  \ \  m\text{-a.e for any} \ f \in L^{p}(E,m).
\end{equation}
We further consider 
\begin{itemize}
\item[(H4)] $(\E,D(\E))$ is conservative.
\end{itemize}

\begin{remark}\label{r;conserv}
Consider the $C_0$-semigroups $(T_t)_{t>0}$, $(\hat{T_t})_{t>0}$ of submarkovian contractions on $L^1(\R^d,m)$. In particular  $(T_t)_{t>0}$ (and also $(\hat{T_t})_{t>0}$) can be defined as semigroups on $L^{\infty}(\R^d,m)$. Then $(\E,D(\E))$ is called conservative, if
\begin{equation}\label{eq;conservative1}
 T_t 1 =1 \ m\text{-a.e. for some (and hence all)} \  t>0 
\end{equation}
Obviously, \eqref{eq;conservative1} holds e.g. if $m(\R^d)<\infty$ and $\|B\|\in L^1(\R^d,m)$. In Section \ref{S5} below we shall present a whole class of examples which do not satisfy 
these two assumptions, but for which \eqref{eq;conservative1}, i.e. (H4) holds.
Clearly \eqref{eq;conservative1} holds, if and only if $m$ is $(\hat{T}_t)$-invariant, that is
\begin{equation}\label{eq;conservative2}
\int \hat{T}_t f \ dm = \int f \ dm \quad \forall f \in L^1(\R^d, m)
\end{equation}
and by \cite[Corollary 2.2]{Stn} \eqref{eq;conservative2} is equivalent to 
\begin{equation}\label{eq;conservative3}
 (1-\hat{L}) \big( C_0^{\infty}(\R^d) \big) \subset L^1(\R^d,m) \ \text{densely}.
\end{equation}
Thus \eqref{eq;conservative3} is equivalent to (H4).
\end{remark}

Following \cite[Proposition 3.8]{AKR}, we obtain:
\begin{prop}\label{prop1.5} 
If (H4) holds (additionally to (H1)-(H3)), then:
\begin{itemize}
\item[(i)] $\lambda R_{\lambda} 1(x) = 1$ for all $x \in E$, $\lambda > 0$.
\item[(ii)] $(P_t)_{t>0}$ is strong Feller on $E$, i.e. $P_t(\mathcal{B}_b(\R^d)) \subset C_b(E)$ for all $t>0$.
\item[(iii)] $P_t 1(x) = 1$ for all $x \in E$, $t>0$.
\end{itemize}
\end{prop} 

By \cite[V. 2.12 (ii)]{MR} (see also \cite[Proposition 1]{Tr5}), it follows that $(\E, D(\E))$ is strictly quasi-regular. Actually, in \cite[section 4.1]{Tr5}, it is shown that this is even true for non-sectorial $B$, i.e. when $\| B \|$ is merely in $L^2_{loc}(\R^d,m)$. In particular, by \cite[V.2.13]{MR} (see also \cite[Theorem 3]{Tr5} for the non-sectorial case) there exists a Hunt process 
$$
\tilde{\bM} = (\tilde{\Omega}, \tilde{\F}, (\tilde{\F})_{t \ge 0}, (\tilde{X}_t)_{t \ge 0}, (\tilde{\P}_x)_{x \in \R^d \cup \{ \Delta \} })
$$ 
with lifetime $\zeta:=\inf\{t\ge 0\,|\,\tilde{X}_t=\Delta\}$ and cemetery $\Delta$ such that $(\E, D(\E))$ is (strictly properly) associated with $\tilde{\bM}$.\\
Consider the strict capacity Cap$_{\E}$ of the non-symmetric Dirichlet form $(\E,D(\E))$ as defined in \cite[V.2.1]{MR} and \cite[Definition 1]{Tr5}, i.e. 
\[
\text{Cap}_{\E} =  \text{cap}_{1,\hat{G}_1\varphi}
\]
for some fixed $\varphi \in L^1(\R^d,m) \cap \mathcal{B}_b(\R^d)$, $0 < \varphi \le 1$. Due to the properties of smooth measures w.r.t. Cap$_{\E}$ in \cite[Section 3]{Tr5} it is possible to consider the work \cite{Tr2} with cap$_{\varphi}$ (as defined in \cite{Tr2}) replaced by Cap$_{\E}$. In particular \cite[Theorem 3.10 and Proposition 4.2]{Tr2} apply w.r.t. the strict capacity Cap$_{\E}$ and therefore the paths of $\tilde{\bM}$ are continuous $\tilde{\P}_x$-a.s. for strictly $\E$-q.e. $x \in \R^d$ on the one-point-compactification $\R^d_{\Delta}$ of $\R^d$ with $\Delta$ as point at infinity.
We may hence assume that 
\begin{equation}\label{contipath}
\tilde{\Omega} = \{\omega = (\omega (t))_{t \ge 0} \in C([0,\infty),\R^d_{\Delta}) \ | \ \omega(t) = \Delta \quad \forall t \ge \zeta(\omega) \}
\end{equation}
and
\[
\tilde{X}_t(\omega) = \omega(t), \quad t \ge 0.
\]
Let Cap be the capacity related to the symmetric Dirichlet form ($\E^0,D(\E^0)$) as defined in \cite[Section 2.1]{FOT}. Then, it holds Cap$(\{\rho=0\})=0$ by \cite[Theorem 2]{fuku85}. 
\begin{lemma}\label{lem2.8}
Let $N \subset \R^d$. Then 
\[
\emph{Cap}(N) = 0  \Rightarrow \emph{Cap}_{\E}(N) = 0.
\]
In particular $\emph{Cap}_{\E}(\{\rho = 0\})=0$.
\end{lemma}
\proof
Let $N \subset \R^d$ be such that Cap$(N)=0$. Then by the definition of Cap
there exist closed sets $F_k \subset \R^d \setminus N$, $k \ge 1$ such that 
\[
\lim _{k \to \infty }\text{Cap}(\R^d \setminus F_k) = 0.
\]
Therefore, we may assume that  Cap$(\R^d \setminus F_k) < \infty$ for any $k \ge 1$. Hence
\[
\mathcal{L}_{\R^d \setminus F_k} : = \{ u \in D(\E^0) \ | \ u \ge 1 \ m\text{-a.e. on} \ \R^d \setminus F_k   \} \neq \emptyset, \quad \forall k \ge 1.
\]
Then by \cite[Lemma 2.1.1.]{FOT} there exists a unique element $e_{\R^d \setminus F_k} \in \mathcal{L}_{\R^d \setminus F_k}$ such that 
\[
\text{Cap}(\R^d \setminus F_k) = \E^0_1(e_{\R^d \setminus F_k}, e_{\R^d \setminus F_k})  \ \ \text{and} \ \ e_{\R^d \setminus F_k} = 1 \ m\text{-a.e on} \ \ \R^d \setminus F_k.
\] 
We denote by $\mathcal{P}$ the family of 1-excessive functions w.r.t. $\E$ in $D(\E)$ and denote by $h_{U}$ the 
(1-) reduced function on an open set $U \subset \R^d$ of a function $h$ in $D(\E)$. Then by (H3) and \cite[III. Proposition 1.5]{MR} for 
$u \le 1, \ u \in \mathcal{P}$
\[
\E_1(u_{\R^d \setminus F_k}, u_{\R^d \setminus F_k} ) \le  \E_1(u_{\R^d \setminus F_k}, e_{\R^d \setminus F_k} )
\le K \  \E_1(u_{\R^d \setminus F_k}, u_{\R^d \setminus F_k} )^{1/2}\E_1(e_{\R^d \setminus F_k}, e_{\R^d \setminus F_k} )^{1/2},
\]
where $K$ is the sector constant.
Therefore,
\[
 \lim_{k \to \infty}  \sup_{\begin{subarray} \ u \le 1, \\ u \in \mathcal{P} \end{subarray}} \E_1(u_{\R^d \setminus F_k}, u_{\R^d \setminus F_k} ) = 0.  
\]
Since for any fixed $\varphi \in L^1(\R^d,m) \cap \mathcal{B}_b(\R^d)$, $0 < \varphi \le 1$
\[
\E_1( u_{\R^d \setminus F_k}, \hat{G}_1 \varphi )  \le K \ \E_1( u_{\R^d \setminus F_k},u_{\R^d \setminus F_k} )^{1/2}   \E_1( \hat{G}_1 \varphi,  \hat{G}_1 \varphi )^{1/2}, 
\]
we have
\[
\text{Cap}_{\E} (N) \le   \lim_{k \to \infty}  \sup_{\begin{subarray} \ u \le 1, \\ u \in \mathcal{P} \end{subarray}} \E_1(u_{\R^d \setminus F_k}, \hat{G}_1 \varphi) = 0.  
\]
\qed

For a Borel set $B \subset \R^d$, we define 
\[
\sigma_B : = \inf \{t > 0 \ | \ \tilde{X}_t \in B  \}, \quad  D_B : = \inf \{t \ge 0 \ | \ \tilde{X}_t \in B  \},
\]
and likewise we define $\sigma_B$, $D_B$ for any other Hunt process. Let 
\[
\tilde{X}_t^{E}(\omega) : = 
\begin{cases}
\tilde{X}_t(\omega) \quad 0 \le t < D_{\R^d \setminus E}(\omega) \\
\Delta \quad t \in [D_{\R^d \setminus E}(\omega), \infty], \ \omega \in \tilde{\Omega}.
\end{cases}
\]
Then $\tilde{\bM}^E : = (\tilde{\Omega}, \tilde{\F}, (\tilde{\F}_t)_{t \ge 0}, (\tilde{X}_t^{E})_{t \ge 0}, (\tilde{\P}_x)_{x \in E \cup \{\Delta\} } )$ is again a Hunt Process by \cite[Theorem A.2.10]{FOT} and its lifetime is $\zeta^E : = \zeta \wedge D_{\R^d \setminus E}$. $\tilde{\bM}^E$ is called the part process of $\tilde{\bM}$ on $E$ and it is associated with the part $(\E^E,D(\E^E))$ of $(\E,D(\E))$ on $E$ (cf. \cite[Theorem 3.5.7]{O13}).
We denote the $L^2(E,m)$-semigroup of $(\E^E,D(\E^E))$ by $(T_t^E)_{t>0}$.

\begin{lemma}\label{nestconti}
Let $(F_k)_{k \ge 1}$ be an increasing sequence of compact subsets of $E$ with $\cup_{k \ge 1} F_k = E$ and such that $F_k \subset   \mathring{F}_{k+1}$, $k \ge 1$(here $\mathring{F}$ denotes the interior of $F$). Then
\[
\tilde{\P}_x (\tilde{\Omega}_0) = 1 \ \text{for strictly} \ \E \text{-q.e.} \ x \in E,
\]
where 
\[
\tilde{\Omega}_0 : = \tilde{\Omega} \cap \{\omega \ | \ \omega(0) \in E \cup \{\Delta\} \ \text{and} \ \lim_{k \to \infty} \sigma_{E \setminus F_k}(\omega) \ge \zeta (\omega) \}.
\]
\end{lemma}
\proof
First note that $\tilde{\P}_x(\zeta = \zeta^E) = 1$ for $m$-a.e. $x \in E$ since Cap$_{\E}(\R^d \setminus E)=0$. By \cite[IV. Theorem 5.1 and Proposition 5.30]{MR} there exists an increasing sequence of compact subsets $(K)_{n \ge 1}$ of $E$ such that 
\[
\tilde{\P}_x(\lim_{n \to \infty} \sigma_{E \setminus K_n} \ge \zeta^E ) = 1\ \text{for} \ m\text{-a.e.} \ x \in E.
\]
The last and previous imply that  
\begin{equation}\label{nest}
\tilde{\P}_x(\lim_{k \to \infty} \sigma_{E \setminus F_k} \ge \zeta ) = 1\ \text{for} \ m\text{-a.e.} \ x \in E
\end{equation}
since $(\mathring{F}_k)_{k \ge 1}$ is an open cover of $K_n$ for every $n \ge 1$. \eqref{contipath} and \eqref{nest} now easily imply  the assertion.

\qed

\begin{thm}\label{existhunt}
There exists a Hunt process
\[
\bM =  (\Omega, \F, (\F_t)_{t \ge 0}, (X_t)_{t \ge 0}, (\P_x)_{x \in E_{\Delta}}   )
\]
with state space $E$, having the transition function $(P_t)_{t \ge 0}$ as transition semigroup. In particular $\bM$ satisfies the absolute continuity condition, because
\[
T_t^E f = P_t f \quad m\text{-a.e.} \ \forall t>0, \ f \in L^2(E,m) \cap \mathcal{B}_b(E).
\]
Moreover $\bM$ has continuous sample paths in the one point compactification $E_{\Delta}$ of $E$ with the cemetery $\Delta$ as point at infinity.
\end{thm}
\proof
Given the transition function $(P_t)_{t \ge 0}$ 
we can construct $\bM$ with continuous sample paths in $E_{\Delta}$ following the line of arguments in \cite{AKR} (see also \cite[Section 2.1.2]{ShTr13a}) using in particular Lemma \ref{nestconti} and our further previous preparations.
As in \cite[Lemma 4.2]{ShTr13a}, we then show that  the (temporally homogeneous) sub-Markovian transition function $(P_t)_{t \ge0}$ on $\big( E, \mathcal{B}(E) \big)$  with transition kernel density $p_t(\cdot, \cdot)$ on $E \times E$ satisfies
\[
T_t^E f = T_t f = P_t f  \quad m\text{-a.e.} 
\]
 for any $t > 0$ and $f \in \mathcal{B}_b(E)$ with compact support (i.e. $|f|dm$ has compact support). Thus the absolute continuity condition is satisfied. 
\qed
\begin{remark}\label{conserve}
If in addition (H4) holds, one can drop $\Delta$ in Theorem \ref{existhunt} and $\bM$ becomes a classical (conservative) diffusion with state space $E$. Indeed, it then holds
\[
\P_x(\zeta = \infty) = 1, \quad \forall x \in E.
\]
\end{remark}

\section{Existence of weak solutions}\label{S3}

\begin{lemma}\label{lem3.1}
Assume (H1)-(H3).
\begin{itemize}
\item[(i)] Let $f \in \bigcup_{s \in [p, \infty)} L^s(E,m)$, $f \ge 0$, then for all $t > 0$, $x \in E$,
\[
\int_0^t P_s f(x) \ ds < \infty,
\]
hence
\[
\int \int_0^t f(X_s) \ ds \ d\P_x < \infty.
\]
\item[(ii)] Let $u \in C_0^{\infty} (\R^d)$, $\lambda > 0$. Then
\[
R_{\lambda} \big( (\lambda - L) u\big)(x) = u(x) \quad \forall x \in E. 
\]
\item[(iii)] Let $u \in C_0^{\infty}(\R^d)$, $t>0$. Then
\[
P_t u(x) - u(x) = \int_0^t P_s(Lu)(x) \ ds \quad \forall x \in E.
\]
\end{itemize}
\end{lemma}
\proof
The proof is the same as the one for \cite[Lemma 5.1]{AKR}. 
\qed

\begin{lemma}\label{lem1.7}
For $u \in C_0^{\infty} (\R^d)$
\[
L u^2 - 2 u \ Lu = \|\nabla u\|^2.
\]
\end{lemma}
\proof
This follows immediately from \eqref{eq1.3}.
\qed

The following result is standard. For the reader's convenience, we include a proof in the Appendix.
\begin{proposition}\label{thm3.2}
Let $u \in C_0^{\infty}(\R^d)$ and
\[
M_t : = \left( u(X_t) - u(X_0) - \int_0^t Lu(X_r) \ dr \right)^2 - \int_0^t \|\nabla u\|^2 (X_r) \ dr, \quad t \ge 0.
\]
Then $(M_t)_{t \ge 0}$ is  an $(\mathcal{F}_t)_{t \ge 0}$-martingale under $\P_x$, $\forall x \in E$.
\end{proposition}

Let $\theta_s : \Omega \rightarrow \Omega$, $s > 0$, be the canonical shift, i.e. $\theta_s(\omega) = \omega(\cdot + s )$, $\omega \in \Omega$.

\begin{lemma}\label{pointwisenest}
Let $(B_k)_{k \ge 1}$ be an increasing sequence of relatively compact open sets in $E$ with $\cup_{k \ge 1} B_k= E$.
Then
for all $x \in E$
\[
\P_x \Big(\lim_{k \rightarrow \infty} \sigma_{E \setminus B_{k}} \ge \zeta \Big)=1.
\] 
\end{lemma}
\proof
Let
\[
\Lambda : = \Big\{ \lim_{k \rightarrow \infty} \sigma_{E \setminus B_{k} }   \ge \zeta \Big\}.
\]
Note that by Lemma \ref{nestconti} for $m$-a.e. $x \in E$
\[
\P_x (\Lambda )=1.
\] 
Then for $x \in E$ and $ s>0$
\begin{eqnarray*}
\P_x(\theta^{-1}_s ( \Lambda)) &= &\EE_x[ 1_{\Lambda} \circ \theta_s]
= \EE_x \Big[  \EE_x[  1_{\Lambda} \circ \theta_s  \ | \ \F_s  ]    \Big] = \EE_x \Big[ \EE_{X_s} [1_{\Lambda}]     \Big]\\
&=& \int_E p_s(x,y)  \ \EE_y[1_{\Lambda}] \ m(dy) + (1-P_s(x, E)) \P_{\Delta}(\Lambda) =1.
\end{eqnarray*}
Let $x \in E$. Define
\[
\Omega_x : = \{ \omega \in \Omega \ | \  t  \mapsto X_t(\omega), \ t \ge0 \ \text{is continuous in } E_{\Delta}  \ \text{and} \ X_0(\omega) = x  \} \ \cap \ \bigcap_{\begin{subarray} \ s>0 \\ s \in S  \end{subarray} }  \ \theta^{-1}_s \circ \Lambda,
\]
where $S$ is a countable dense set in $(0, \infty)$.
Fix $\omega \in \Omega_x$. By the continuity of $X_t(\omega)$ there is $s^{\prime} \in S $ such that $X_t(\omega) \in B_{\bar{k}}$, $t \in [0,s^{\prime}]$, for some $\bar{k} \in \N$.   This implies
\[
\sigma_{E \setminus B_k}(\omega) = s^{\prime} + \sigma_{E \setminus B_k} (\theta_{s^{\prime}}(\omega) )
\]
for $k \ge \bar{k}$ and since $\zeta(\omega) \ge s^{\prime}$, we get
\[
\zeta(\omega) = s^{\prime} + \zeta (\theta_{s^{\prime}}(\omega)).
\]
Putting all together and noting that $\theta_{s^{\prime}}(\omega) \in \Lambda$, we obtain
\[
\lim_{k \to \infty} \sigma_{E \setminus B_k}(\omega) = \lim_{k \to \infty} \sigma_{E \setminus B_k}( \theta_{s^{\prime}}(\omega)) + s^{\prime} 
\ge \zeta(\theta_{s^{\prime}}(\omega)) + s^{\prime}  = \zeta(\omega).
\]
Hence $\Omega_x\subset\Lambda$. Since $\P_x( \Omega_x)=1$, the assertion follows. 
\qed

\begin{remark}
For an alternative proof of Lemma \ref{pointwisenest}, which does not require the absolute continuity condition, we refer to Lemma \ref{l;localize} in Section \ref{Appendix}.
\end{remark}

\begin{thm}\label{t3.6}
Under (H1)-(H3) after enlarging the stochastic basis $(\Omega, \F, (\F_t)_{t\ge 0},\P_x )$ appropriately for every $x\in E$, the process $\bM$ satisfies 
\begin{eqnarray}\label{weaksolution}
X_t = x + W_t +  \int_0^t  \left( \frac{\nabla \rho}{2 \rho} + B \right) (X_s) \ ds, \quad t < \zeta
\end{eqnarray}
$\P_x$-a.s. for all $x \in E$ where $W$ is a standard $d$-dimensional $(\F_t)$-Brownian motion on $E$. If additionally (H4) holds, 
then we do not need to enlarge the stochastic basis and $\zeta$ can be replaced by $\infty$ (cf. Remark \ref{conserve}).
\end{thm}
\proof
Let $u_i \in C_0^{\infty}(E)$, $i=1,\dots,d$, and 
\[
M_t^{u_i} : = u_i(X_t) - u_i(X_0) - \int_0^t Lu_i(X_s) \ ds, \quad 1\le i \le d, \ t \ge 0.
\]
For $x \in E$, $(M^{u_i}_{t})_{t \ge 0}$ is a continuous $(\mathcal{F}_t)_{t \ge 0}$-martingale under $\P_x$. Note that by Proposition \ref{thm3.2} and polarization, the quadratic covariation processes satisfy
\[
\langle M^{u_i}, M^{u_j} \rangle_t =  \int_0^t  \langle \nabla u_i, \nabla u_j \rangle (X_s) \ ds, \quad 1 \le i,j  \le d, \ t \ge 0.
\]
Suppose $\zeta<\infty$. Then there is an enlargement $(\bar{\Omega}, \bar{\F}, \bar{\P}_x )$ (since $\langle \nabla u_i, \nabla u_j \rangle 1_E$ is degenerate on $E_{\Delta}$) of the underlying probability space $(\Omega, \F, \P_x )$, 
a d-dimensional Brownian motion $(W_t)_{t \ge 0} = (W_t^1,\dots, W_t^d)_{t \ge 0}$ on $(\bar{\Omega}, \bar{\F}, \bar{\P}_x )$ and a $d \times d$ matrix $\sigma = (\sigma_{ij})_{1 \le i,j \le d}$ such that
\[
M_t^{u_i} = \sum_{k=1}^{d}  \int_0^t \ \sigma_{ik} (X_s) \ dW_s^k, \quad 1 \le i \le d, \  t\ge 0.
\]
and $\langle \nabla u_i, \nabla u_j \rangle = \sum_{k=1}^{d} \sigma_{ik} \sigma_{jk}$ (cf. \cite[Section 3.4.A., 4.2 Theorem]{KS}).
The identification of $X$ up to $\zeta$ is now obtained by using Lemma \ref{pointwisenest} with an appropriate localizing sequence as in Lemma \ref{nestconti}
for which the coordinate projections on $E$ coincide locally with $C_0^{\infty}(E)$-functions and noting that $W_t^i=\int_0^t 1_E(X_s)dW_s^i$ on $\{t<\zeta\}$. 
If $\zeta=\infty$, using the same localization, we obtain that $\langle M^{v_i} \rangle_t =  \int_0^t  1_E(X_s) \ ds=t$ for $t<\infty$, 
where $v_i$ is the i-th coordinate projection. Thus $M^{v_i}$  is a Brownian motion by L\'evy's characterization and we do not need an enlargement of the stochastic basis. The localization of the drift part is trivial.
\qed

\section{Pathwise uniqueness and strong solutions}\label{S4}
We first recall that by \cite[Theorem 2.1]{KR} under the conditions (H1), (H2) ((H3) is not needed), for every stochastic basis and given Brownian motion $(W_t)_{t\geq0}$ there 
exists a strong solution to \eqref{weaksolution} which is pathwise unique
among all solutions satisfying
\begin{equation}\label{eq4.1}
 \int_0^t \left\| \left( \frac{\nabla \rho}{2 \rho} + B \right)(X_s)\right\|^2 \mathrm{d} s <\infty \quad \mathbb{P}_x\textnormal{-a.s. on } \{t<\zeta\}\;.
\end{equation}
In addition, one has pathwise uniqueness and weak uniqueness in this class. 

In the situation of Theorem \ref{t3.6} it follows, however immediately from Lemma \ref{pointwisenest} that \eqref{eq4.1} holds for the solution there.  Indeed, by Lemma \ref{pointwisenest}, (\ref{eq4.1}) holds with $\sigma_{E\setminus B_k}$ 
for all $k\in \N$. But the latter together with (H1) clearly implies that (\ref{eq4.1}) holds $\P_x$-a.s. for all $x\in S$ for some $S\in {\cal B}(E)$ with $m(E\setminus S)=0$  (by Lemma \ref{lem2.8} the set $S$ can be chosen such that even $\text{Cap}_{\E}(E\setminus S)=0$). 
So, \cite[Theorem 2.1]{KR}, in particular, implies that the law of 
$\tilde \P_x$ of the strong solution from that theorem coincides with $\P_x$ for all $x\in S$. But then $\tilde \P_x= \P_x$ 
for all $x\in E$, because of the strong Feller property of our Markov process given by $(\P_x)_{x\in E}$ and of the 
one from \cite[Theorem 2.1]{KR}, i.e. $\tilde \P_x$, $x\in E$, since $S$ is dense in $E$. In particular, (\ref{eq4.1}) 
holds for all $x\in E$. Hence we obtain the following:
\begin{thm}\label{t4.1}
 Assume (H1)-(H3). For every $x \in E$ the solution in Theorem \ref{t3.6} is strong, pathwise and weak unique. In particular, it is adapted to the filtration $(\mathcal{F}_t^W)_{t\geq0}$ generated by the Brownian motion $(W_t)_{t\geq0}$
 in \eqref{weaksolution}. 
\end{thm}
\begin{remark}\label{r4.2}
(i) By Theorem \ref{t3.6} and \ref{t4.1} we have thus shown that (the closure of) \eqref{df} 
	is the Dirichlet form associated to the Markov processes given by the laws of the (strong) solutions to \eqref{weaksolution}.
  Hence we can use the theory of Dirichlet forms to show further properties of the solutions. \\ 
	(ii) In \cite{KR} also a new non-explosion criterion was proved (hence one obtains (H4)), assuming that $\frac{\nabla \rho}{2 \rho}+B$ is the (weak) gradient of a function 
	$\psi$ which is a kind of Lyapunov function for \eqref{weaksolution}.
  The theory of Dirichlet forms provides a number of analytic non-explosion, i.e. conservativeness criteria 
	(hence implying (H4)) which are completely different from the usual ones for SDEs and which are checkable in many cases.
  As stressed in (i) such criteria can now be applied to \eqref{weaksolution}. Even the simple already mentioned case, where $m(\R^d)<\infty$ and $\|B\| \in L^1(\R^d,m)$ which entails (H4), appears to be a new non-explosion condition for \eqref{weaksolution}. Further explicit examples where \eqref{weaksolution} has a non-explosive unique strong solution are given in Section \ref{S5} below.

\end{remark}

\section{Applications to Muckenhoupt $A_{\beta}$-weights}\label{S5}
In this section we present a class of examples of $\rho$ and $B$ satisfying our assumptions (H3) and (H4). Throughout, we assume (H1) and (H2) to hold.
\begin{lemma}\label{t;sector}
Suppose 
\begin{itemize}
\item[(i)] For $r>0$
\[
\left(  \int_{B_r(0)} | u |^{\frac{2N}{N-2}} \ \rho \ dx \right) ^{\frac{N-2}{2N}} \le c_r \left( \int_{B_{2r}(0)}  \left( \|\nabla u \|^2 + u^2 \right) \ \rho \ dx    \right)^{1/2}, \quad \forall u \in C_0^{\infty}(\R^d),
\]
where $c_r$ is some constant, $N>2$ and
\item[(ii)] $\|B\| \in L^{N}_{loc}(\R^d,m)  \cap L^{\infty}(K^c,m)$ for some compact $K \subset \R^d$. 
\end{itemize}
 Then
\[
\left| \int_{\R^d} \langle B,\nabla u \rangle \ v \ \rho \ dx \right| \le c_{B,K} \ \E_1^0 (u,u)^{1/2} \ \E_1^0 (v,v)^{1/2}, \quad \forall u,v \in C_0^{\infty}(\R^d),
\]
where $c_{B,K}$ is some constant, i.e. (H3) holds.
\end{lemma}
\proof
For $r_0 >0$ such that $K \subset B_{r_0}(0)$
\begin{eqnarray*}
&& \left| \int_{\R^d} \langle B, \nabla u \rangle \ v \  \rho \ dx \right| \le \left( \int_{\R^d} \|B \|^2 \ v^2 \ \rho \ dx   \right)^{1/2} \ \left( \int_{\R^d}  \| \nabla u \|^2  \rho \ dx   \right)^{1/2} \\
&\le&  \left( \int_{B_{r_0} (0) } \|B \|^2 \ v^2 \ \rho \ dx  +  \int_{B_{r_0}(0)^c } \|B \|^2 \ v^2 \ \rho \ dx           \right)^{1/2} \ \E_1 (u,u)^{1/2}\\
&\le& \left(   \left( \int_{B_{r_0} (0) } \|B \|^2 \ v^2 \ \rho \ dx \right)^{1/2} +   \| B \|_{\infty, K^c } \| v \|_{L^2(\R^d, m)}         \right) \  \E_1 (u,u)^{1/2} \\
&\le& \left(   \left( \int_{B_{r_0} (0) } \|B \|^N  \ \rho \ dx \right)^{1/N}  \left( \int_{B_{r_0} (0) } v^{ \frac{2 N}{N-2}}  \ \rho \ dx \right)^{  \frac{N-2}{2N}  }  +   \| B \|_{\infty, K^c} \| v \|_{L^2(\R^d, m)}         \right) \  \E_1 (u,u)^{1/2} \\
&\le& c_{B,K} \ \E_1^0 (u,u)^{1/2} \ \E_1^0 (v,v)^{1/2}.
\end{eqnarray*}
The last inequality follows from assumption (i) and $\| \cdot \|_{\infty,K^c}$ denotes the $L^{\infty}(\R^d,m)$-norm on $K^c$.
\qed

\begin{lemma}\label{l;sobolmu}
Let  $\rho$ be a Muckenhoupt $A_{\beta}$-weight, $1  \le \beta \le 2$. Then for $x \in \R^d$, $r>0$, $N>2$
\[
\left(  \int_{B_r(x)} | u |^{\frac{2N}{N-2}} \ dm \right) ^{\frac{N-2}{2N}} \le C_{x,r} \left( \int_{B_{2r}(x)}  \left( \|\nabla u \|^2 + u^2 \right) \ dm   \right)^{1/2}, \quad \forall u \in C^{\infty}(\R^d), 
\]
where $C_{x,r}$ is some constant and  $N \ge \beta d + \log_2 A$, $A$ is the $A_{\beta}$ constant of $\rho$. 
\proof
By the doubling property of $A_{\beta}$-weights (cf. \cite[Proposition 1.2.7]{Tu} ), 
\begin{equation}\label{eq;doubling}
m (B_{2r}(x)) \le A \ 2^{\beta d} \ m (B_{r}(x)).
\end{equation}
Note that $A_{\beta} \subset A_2$ if $1 \le \beta \le 2$. Then
by \cite[Theorem (1.5)]{FKS}  the scaled Poincar\'{e} inequality holds true, i.e.
for $x \in \R^d$, $r>0$
\[
\int_{B_r(x)} | u - u_{x,r} |^2  \ dm   \le c r^2 \int_{B_{r}(x)}   \|\nabla u \|^2  \ dm, \quad \forall u \in C^{\infty}(\R^d),    
\]
where $u_{x,r} = \frac{1}{m(B_r(x))}\int_{B_r(x)} u \ dm$ and $c$ is some constant.
Consequently, \cite[Theorem 2.1]{Sc},  the doubling property, and the scaled Poincar\'{e} inequality imply the Sobolev inequality, i.e.
for $x \in \R^d$, $r>0$, $N >2$
\[
\left(  \int_{B_r(x)} | u |^{\frac{2N}{N-2}} \ dm \right) ^{\frac{N-2}{2N}} \le c_{x,r} \left( \int_{B_{r}(x)}  \left( \|\nabla u \|^2 + u^2 \right) \ dm    \right)^{1/2}, \quad \forall u \in C_0^{\infty}(B_r(x)), 
\]
where $c_{x,r}$ is some constant and $N \ge \beta d + \log_2 A$. Then using a cutoff function like for instance $g_r(y) : = \frac{1}{r} (2r - \| x-y\|)^+$, we see that for $x \in \R^d$, $r>0$  
\[
\left(  \int_{B_r(x)} | u |^{\frac{2N}{N-2}} \ dm \right) ^{\frac{N-2}{2N}} \le C_{x,r} \left( \int_{B_{2r}(x)}  \left( \|\nabla u \|^2 + u^2 \right) \ dm    \right)^{1/2}, \quad \forall u \in C^{\infty}(\R^d), 
\]
where $C_{x,r}$ is some constant and $N > 2$ as well as $N \ge \beta d + \log_2 A$. 
\qed
\end{lemma}

\begin{lemma}\label{l;bsector}
Let  $\rho$ be a Muckenhoupt $A_{\beta}$ weight, $1  \le \beta \le 2$, $N>2$ and $\|B\| \in L^{N}_{loc}(\R^d,m)  \cap L^{\infty}(K^c,m)$ for some compact $K \subset \R^d$, $N \ge \beta d + \log_2 A$,  where $A$ is the $A_{\beta}$ constant of $\rho$. 
 Then
\[
\left| \int_{\R^d} \langle B,\nabla u \rangle \ v \ \rho \ dx \right| \le c_{B,K} \ \E_1^0 (u,u)^{1/2} \ \E_1^0 (v,v)^{1/2}, \quad \forall u,v \in C_0^{\infty}(\R^d),
\]
where $c_{B,K}$ is some constant, i.e. (H3) holds.
\end{lemma}
\proof
This follows from Lemma \ref{t;sector} and Lemma \ref{l;sobolmu}.
\qed

\begin{lemma}\label{l;muconserv}
It holds 
\[
(1- \hat{L})(C_0^{\infty}(\R^d)) \subset L^1 (\R^d, m)  \ \ \text{densely}.
\]
In particular (H4) holds (cf. Remark \ref{r;conserv}).
\end{lemma}
\proof
Let $h \in L^{\infty}(\R^d,m)$ be arbitrary. We have to show that 
\begin{equation}\label{eq;hb}
\int (1-  \hat{L}) f \cdot h \ dm = 0 \quad \forall f \in C_0^{\infty}(\R^d)
\end{equation}
implies $h=0$.\\
By \cite[Theorem 2.1]{Stn} it follows from \eqref{eq;hb} that $h \in D(\E^0)_{loc}: = \{ u \ | \ u \cdot \chi \in D(\E^0) \ \forall  \chi   \in C_0^{\infty}(\R^d)   \}$ and 
\begin{equation}\label{eq;diriequa}
\E_1^0 (u,h) = - \int \langle B, \nabla u \rangle h \ dm \quad \forall u \in D(\E^0)_0
\end{equation}
where $D(\E^0)_0 := \{ u \in D(\E^0) \ | \ \  \text{supp} (|u| dm) \ \  \text{is compact}      \}$.
Define
\begin{eqnarray*}
\varv (r) :&=& m(B_r(0)),  \quad r>0 \\
a_n : &=& \int_n^{2n} \frac{s}{\log (\varv(s)) } \ ds, \quad n \ge 1\\
\psi_n(r) : &=& 1_{[0,n]}(r) - \frac{1}{a_n} \int_n^r \frac{s}{\log(\varv(s))} \ ds \cdot 1_{[n,2n]} (r)\\
u_n(x) : &=& \psi_n(\|x\|).
\end{eqnarray*}
Then $u_n \in D(\E^0)_0$ and
\begin{eqnarray}
\nabla u_n (x) &= &- \frac{1}{a_n} \frac{x}{\log (\varv (\|x\|))} \cdot 1_{[n,2n]} (\|x\|)   \label{eq;grun} \\
a_n &\ge& \int_n^{2n} \frac{n}{\log(\varv(2n))} \ ds = \frac{n^2}{\log(\varv(2n))}  \ge \frac{n^2}{\log( A 2 ^{\beta d}  \ \varv(n))} .\label{eq;anine} 
\end{eqnarray}
The last inequality follows from \eqref{eq;doubling}.  Taking sufficiently large $n$ such that  $\log( A 2 ^{\beta d}) \le \log(  \varv(n)), $
\eqref{eq;grun} and \eqref{eq;anine}  imply
\begin{equation}\label{eq;grun2}
\|\nabla u_n(x)\| \le  \frac{\log( A 2 ^{\beta d}  \ \varv(n))}{n^2} \ \frac{2n}{\log(\varv(n))} \cdot 1_{[n,2n]} (\|x\|) \le \frac{4}{n} \cdot 1_{[n,2n]}(\|x\|).
\end{equation}
Then
\begin{eqnarray*}
\phi(n) : &=& \int_{B_n(0)} h^2 \ dm \le \int_{B_{2n}(0)} h^2 u_n^2 \ dm = \int_{B_{2n}(0)} (h u_n^2)  \cdot h \ dm \\
&=& - \int_{B_{2n}(0)}   \langle \nabla ( h u_n^2), \nabla h \rangle \ dm - \int_{B_{2n}(0)}    \langle B, \nabla (h u_n^2) \rangle   h  \ dm
\end{eqnarray*}
Since $h u_n^2 \in D(\E^0)_0$, the last equality follows from \eqref{eq;diriequa}. The last term is equal to
\[
- \int_{B_{2n}(0)}    \langle B, \nabla (u_n \cdot (h u_n)) \rangle   h  \ dm = 
- \int_{B_{2n}(0)}    \langle B, \nabla (u_n ) \rangle   h^2 u_n  \ dm
- \int_{B_{2n}(0)}    \langle B, \nabla  (h u_n)  \rangle   h u_n  \ dm.
\]
Since $h u_n \in D(\E^0)_0$, the second term is zero by (H3). Therefore 
\begin{eqnarray*}
\phi(n) &=& - \int_{B_{2n}(0)}   \langle \nabla ( h u_n^2),  \nabla h \rangle  \ dm
- \int_{B_{2n}(0)}    \langle B, \nabla (u_n ) \rangle   h^2 u_n  \ dm\\
&=&   - \int_{B_{2n}(0)}   \big( u_n^2 \|\nabla h\|^2 + 2 u_n \langle \nabla h, \nabla u_n \rangle h    \big) \ dm
- \int_{B_{2n}(0)}    \langle B, \nabla (u_n ) \rangle   h^2 u_n  \ dm\\
&\le&  \int_{B_{2n}(0)} \| \nabla u_n \|^2 h^2 \ dm  - \int_{B_{2n}(0)}    \langle B, \nabla (u_n ) \rangle   h^2 u_n  \ dm.
\end{eqnarray*}
Taking $n \ge 4$ so large that $K \subset B_n(0)$ and that \eqref{eq;grun2} holds
\begin{eqnarray*}
\phi(n) &\le&   \left( \frac{4}{n} \right)^2 \int_{B_{2n}(0) \setminus B_n(0)} h^2 \ dm + \frac{4}{n} \| B \|_{\infty, K^c} \int_{B_{2n}(0) \setminus B_n(0)} h^2 \ dm \\
&\le& \frac{4}{n}  (\| B \|_{\infty, K^c} +1) \int_{B_{2n}(0) } h^2 \ dm = \frac{4}{n} (\| B \|_{\infty, K^c} +1)  \phi(2n).
\end{eqnarray*}
Set $C:=  4 (\| B \|_{\infty, K^c} +1)$. Thus by iteration of the last inequality and \eqref{eq;doubling}, we obtain for any $k \ge1$
\[
\phi(n) \le \frac{C^k}{n^k 2^{ \frac{k(k+1)}{2}    }  } \  \phi(2^kn)  \le 
\frac{C^k}{n^k 2^{ \frac{k(k+1)}{2}    }  } \  \|h\|_{\infty}^2  \ \varv(2^k n)
\le \frac{C^k}{n^k 2^{ \frac{k(k+1)}{2}    }  } \|h\|_{\infty}^2   \   (A 2^{\beta d})^k \  \varv( n).
\]
Note that $\varv(n) \le c n^{\alpha}$ for some $\alpha > 0$, where $c>0$ is some constant. Now choose $k > \alpha$ then $\phi(n) \to 0$ as $n \to \infty$, hence $h=0$.
\qed

Lemma \ref{l;bsector} and Lemma \ref{l;muconserv} imply the final theorem.

\begin{thm}\label{Muckenhouptmain}
Let $\rho$ and $B$ satisfy the assumptions (H1) and (H2) and the assumptions of Lemma \ref{l;bsector}. Then (H1)-(H4) hold. Consequently, Theorems \ref{t3.6} and \ref{t4.1} apply with $\zeta=\infty$.
\end{thm}

\section{Appendix}\label{Appendix}

\proof (of Proposition \ref{thm3.2})
By Lemma \ref{lem3.1} and the Markov property
\[
u(X_t) - u(X_0) - \int_0^t Lu(X_r) \ dr, \quad u \in C_0^{\infty}(\R^d), \  t \ge 0
\]
is a square integrable $(\mathcal{F}_t)_{t \ge 0}$-martingale under $\P_x$ for all $x \in E$. 
Fix $x \in E$,  $u \in C_0^{\infty}(\R^d)$, and set
\[
M_t : = \left( u(X_t) - u(X_0) - \int_0^t Lu(X_r) \ dr \right)^2 - \int_0^t \|\nabla u\|^2 (X_r) \ dr, \quad t \ge 0.
\]
Then since $u \in D(L_p)$ (cf. \eqref{eq1.3}), it follows by Lemma \ref{lem3.1} that  $(M_t)_{t \ge 0}$ and all integrands below are integrable w.r.t. $\P_x$.
Using Lemma \ref{lem1.7} we get for $s \in [0,t)$
\begin{eqnarray*}
&& M_t - M_s\\
 &= &\left( u(X_t) - u(X_0) - \int_0^t Lu(X_r) \ dr  + u(X_s) - u(X_0) - \int_0^s Lu(X_r) \ dr \right) \\
&& \times \left( u(X_t) - u(X_s) - \int_s^t Lu(X_r) \ dr\right) - \int_s^t ( L u^2 - 2u \ Lu) (X_r) \ dr\\
 &= &\left( u(X_t) + u(X_s) - 2 u(X_0) -2 \int_0^s Lu(X_r) \ dr  - \int_s^t Lu(X_r) \ dr \right) \\
&& \times \left( u(X_t) - u(X_s) - \int_s^t Lu(X_r) \ dr\right) - \int_s^t ( L u^2 - 2u \ Lu) (X_r) \ dr\\
&=& u^2(X_t) - u^2(X_s) -2u  (X_0) \Big(u(X_t) - u(X_s)\Big) \\
&& -2\Big(u(X_t) - u(X_s)\Big) \int_0^s Lu (X_r) \ dr - \Big(u(X_t) - u(X_s)\Big) \int_0^{t-s} Lu(X_{r+s}) \ dr\\
&& -\Big(u(X_t) + u(X_s)\Big) \int_0^{t-s} Lu(X_{r+s}) \ dr + 2 u(X_0) \int_0^{t-s} Lu(X_{r+s}) \ dr\\
&&+ 2 \int_0^s Lu (X_r) \ dr \int_0^{t-s} Lu(X_{r+s}) \ dr + \left( \int_0^{t-s} Lu(X_{r+s}) \ dr \right)^2\\
&&-\int_s^t \Big( Lu^2 - 2u \ Lu \Big) (X_r) \ dr.
\end{eqnarray*}
Taking conditional expectation, it follows $\P_x$-a.s.
\begin{eqnarray*}
&&\EE_x [M_t - M_s \ | \ \mathcal{F}_s] =  P_{t-s}u^2 (X_s) - u^2 (X_s) \\
&& - 2u(x) \Big(P_{t-s} u(X_s) - u(X_s) \Big) - 2 \Big(P_{t-s} u(X_s) - u(X_s)  \Big) \int_0^s L u(X_r) \ dr\\
&& -2 \EE_x \Big[u(X_t) \int_0^{t-s} Lu(X_{r+s}) \ dr \ | \ \mathcal{F}_s\Big] + 2 u(x) \int_0^{t-s} P_r(L u) (X_s) \ dr\\
&& + 2 \int_0^s Lu (X_r) \ dr \int_0^{t-s} P_r (Lu) (X_s) \ dr + \EE_x \Big[ \left(\int_0^{t-s} Lu(X_{r+s}) \ dr \right)^2 \ | \ \mathcal{F}_s \Big]\\ 
&&- \int_0^{t-s} P_r \Big(L u^2 - 2u \ Lu  \big) (X_s) \ dr.
\end{eqnarray*}
Using  Lemma \ref{lem3.1}(iii) this simplifies to 
\begin{eqnarray*}
&&\EE_x [M_t - M_s \ | \ \mathcal{F}_s] =  -2 \EE_x \Big[u(X_t) \int_0^{t-s} Lu(X_{r+s}) \ dr \ | \ \mathcal{F}_s\Big] \\
&& + \EE_x \Big[ \left(\int_0^{t-s} Lu(X_{r+s}) \ dr \right)^2 \ | \ \mathcal{F}_s \Big]  + 2\int_0^{t-s} P_r \Big( u \ Lu \Big)  (X_s) \ dr.
\end{eqnarray*}
Note that the first term of the right hand side satisfies
\[
-2 \EE_x \Big[u(X_t) \int_0^{t-s} Lu(X_{r+s}) \ dr \ | \ \mathcal{F}_s\Big]  = -2 \int_0^{t-s} P_r \Big(Lu \ P_{t-s-r} u  \Big)(X_r) \ dr 
\]
and the second term  satisfies
\begin{eqnarray*}
&& \EE_x \Big[ \left(\int_0^{t-s} Lu(X_{r+s}) \ dr \right)^2 \ | \ \mathcal{F}_s \Big]  = 2 \int_0^{t-s} \int_0^{r^{\prime}} \EE_{X_s} \Big[Lu(X_r) Lu(X_{r^{\prime}}) \Big] \ dr \ dr^{\prime}\\
&=& 2 \int_0^{t-s} \int_0^{r^{\prime}} P_r \Big(Lu P_{r^{\prime} -r } (Lu) \Big)(X_s)  \ dr \ dr^{\prime}\\
&=& 2 \int_0^{t-s} P_r \Big(Lu \big(P_{t-s-r} u - u   \big) \Big)(X_s) \ dr
\end{eqnarray*}
by Fubini's theorem. Therefore $\EE_x [M_t - M_s \ | \ \mathcal{F}_s] = 0$  $\P_x$-a.s. and the assertion follows.
\qed

We present here an alternative proof of Lemma \ref{pointwisenest}, which does not require the absolute continuity condition.
\begin{lemma}\label{l;localize}
Let $(B_k)_{k \ge 1}$ be an increasing sequence of relatively compact open sets  in $E$ with $\cup_{k \ge 1} B_k= E$. Then
for all $x \in E$
\[
\P_x \Big(\lim_{k \rightarrow \infty} \sigma_{E \setminus B_{k}} \ge \zeta \Big)=1.
\] 
\end{lemma}
\proof
Let $(B_k)_{k \ge 1}$ be an increasing sequence of relatively compact open sets in $E$ with $\cup_{k \ge 1} B_k= E$ and $\sigma : = \lim_{k \rightarrow \infty} \sigma_{E \setminus B_{k}}$. By quasi-left-continuity of $\bM$ 
\begin{equation}\label{qlc}
\P_x \Big(\lim_{k \rightarrow \infty} X_{\sigma_{E \setminus B_{k}}} =X_{\sigma}, \ \sigma < \infty \Big) = \P_x(\sigma < \infty), \quad \forall x \in E.
\end{equation}
Using \cite[Lemma A.2.7]{FOT}, it follows for any $k \ge 1$ that  
\[
\P_x\Big( X_{\sigma_{E \setminus B_k}} \in (E \setminus B_k) \cup \{\Delta\}, \ \sigma_{E \setminus B_k} < \infty \Big) = \P_x(\sigma_{E \setminus B_k} < \infty), \quad \forall x \in E,
\] 
hence
\begin{equation}\label{huntpro}
\P_x \Big( X_{\sigma_{E \setminus B_k}} \in (E \setminus B_k) \cup \{\Delta\}, \ \sigma < \infty \Big) = \P_x(\sigma < \infty), \quad \forall x \in E.
\end{equation}
From \eqref{qlc} and \eqref{huntpro}
\[
\P_x \Big(\lim_{k \rightarrow \infty} X_{\sigma_{E \setminus B_{k}}} =X_{\sigma}, \  X_{\sigma_{E \setminus B_k}} \in (E \setminus B_k) \cup \{\Delta\},\ \forall k \ge 1 ,\ \sigma < \infty \Big) = \P_x(\sigma < \infty), \quad \forall x \in E.
\]
Let 
\[
A : =   \Big\{ \lim_{k \rightarrow \infty} X_{\sigma_{E \setminus B_{k}}} =X_{\sigma}, \  X_{\sigma_{E \setminus B_k}} \in (E \setminus B_k) \cup \{\Delta\}, \  \forall k \ge 1, \  \sigma < \infty \Big\}, \quad B : = \Big\{X_{\sigma} \in \{\Delta\} \Big\}.
\]
Suppose, to show $A \subset B$, that $\omega \in A$ but $\omega \notin B$, i.e. there exists  $x \in E$ such that $X_{\sigma(\omega)} (\omega) = x $ with $\omega \in A$ .  Since $E$ is open in $\R^d$, we can find a ball $B_{\varepsilon}(x)$, $\varepsilon > 0$ such that the closure $\overline{B_{\varepsilon}(x)} \subset E$. Since $(B_k)_{k \ge 1}$ is an open cover of $\overline{B_{\varepsilon}(x)}$ and  increasing, we can find $k^{\star} \in \N$ such that $B_k \supset \overline{B_{\varepsilon}(x)}$ for all $k \ge k^{\star}$. Since $\omega \in A$, this implies that $X_{\sigma_{E \setminus B_{k}}(\omega) }(\omega) \notin B_{\varepsilon}(x)$, $k \ge k^{\star}$ and so  $ \lim_{k \rightarrow \infty} X_{\sigma_{E \setminus B_{k}}(\omega)}(\omega)  \notin B_{\varepsilon}(x)$, which draws a contradiction.
Hence
\[
\P_x \Big( X_{\sigma} \in \{\Delta\}, \ \sigma < \infty \Big) = \P_x(\sigma < \infty ), \quad \forall x \in E,
\]
and so
\[
\P_x \Big( \sigma \ge \zeta, \ \sigma < \infty \Big) = \P_x(\sigma < \infty), \quad \forall x \in E.
\]
Clearly
\[
\P_x \Big( \sigma \ge \zeta, \ \sigma = \infty \Big) = \P_x(\sigma = \infty), \quad \forall x \in E,
\]
thus 
\[
\P_x( \sigma \ge \zeta) = 1, \quad \forall x \in E.
\]
\qed

Michael R\"ockner\\ 
Fakult\"at f\"ur Mathematik\\
Universit\"at Bielefeld\\
Universit\"atsstrasse 25\\
33615 Bielefeld, Germany, \\
E-mail: roeckner@math.uni-bielefeld.de\\ \\
Jiyong Shin\\
School of Mathematics\\
Korea Institute for Advanced Study\\
85 Hoegiro Dongdaemun-gu, \\
Seoul 02445, South Korea, \\
E-mail: yonshin2@kias.re.kr\\ \\
Gerald Trutnau\\
Department of Mathematical Sciences and \\
Research Institute of Mathematics of Seoul National University,\\
1, Gwanak-Ro, Gwanak-Gu \\
Seoul 08826, South Korea,  \\
E-mail: trutnau@snu.ac.kr
\end{document}